\documentclass[final,letterpaper,oneside,10pt]{article}

\author{David~S.~Knight \\  \texttt{davids.knight42@gmail.com}}
\date{}
\title{On Calculating Square Roots in $GF(p)$}

\begin{document}
\maketitle

\begin {abstract} 
This article presents a new method for calculating square roots in $GF(p)$
by exponentiating in $GF(p^3)$ or equivalently modulo irreducible cubic polynomials.  This algorithm
is in some ways similar to the Cipolla-Lehmer algorithm which is based on exponentiating in $GF(p^2)$.
Another less well known square root algorithm based on quadratic sums is also given.  In addition to this,
several conjectures about the output of this $GF(p^3)$ square root algorithm are mentioned.

\end{abstract}

\section*{1. Introduction}
The two most well-known algorithms for calculating square roots in $GF(p)$ are the Cipolla-Lehmer and the 
Tonelli-Shanks algorithms, both of which are described in [9].  Some variations of the Tonelli-Shanks algorithm
are described in [6].  The Cipolla-Lehmer algorithm is asymptotically the fastest and runs in $O(M(p)$ $log$  $p)$ time
where $M(p)$ is the amount of time it takes to calculate one multiplication modulo $p$.  The Tonelli-Shanks is slower 
and runs in $O(M(p)~v~log~p)$ time where $v$ is the greatest integer such that $2^v$ divides $p-1$.
However, if $v$ is small, the Tonelli-Shanks algorithm is generally faster since it is based on exponentiation in $GF(p)$
whereas the Cipolla-Lehmer exponentiates in $GF(p^2)$ which is less efficient.

Most of the other algorithms for calculating square roots in $GF(p)$ are in some way based on one of these two 
algorithms.  A few algorithms, however are not in any way related to either of these two algorithms.  Two examples 
are the algorithms of Schoof [10] or of Sze [11], both of which use elliptic curves in different ways to calculate square
roots in $GF(p)$.  Another example is the algorithm mentioned in [1] which uses quadratic sums in $GF(p)$
in order to calculate square roots in certain cases.

A new algorithm presented in this paper uses exponentiation in $GF(p^3)$ in order to calculate square roots in $GF(p)$.
While this algorithm is significantly different from previously known methods, it is more closely related to the Cipolla-Lehmer than
it is to the Tonelli-Shanks method.  Part of this algorithm depends on calculating cube roots.  A standard method for calculating 
cube roots based on Tonelli-Shanks  results in an  $O(M(p)(t+1)~log~p)$ algorithm where $t$ is the greatest
integer such that $3^t$ divides $p-1$.  However more efficient methods for calculating cube roots exist, for example 
see [2] or [9].  This means that this $GF(p^3)$ square root algorithm can actually be implemented to run in $O(M(p)~log~p)$ time.
This square root algorithm for the case $p \equiv 5$ $(mod$ $6)$ is presented in Section 4 as Algorithm 2.

In this paper two previously known methods for calculating square roots in $GF(p)$ will be mentioned.  First in Section 2 the 
quadratic sum method is given.  Then in Section 3 the Cipolla-Lehmer algorithm is mentioned.  In Section 4, a new algorithm
based on exponentiation in $GF(p^3)$ is presented and  in Section 5, some conjectures that are related to this algorithm are given.
\\

\section*{2. The Quadratic Sum Method}
In [1] Bach and Huber describe a method for calculating square roots in $GF(p)$ in certain cases based on quadratic sums, which can be
considered to be a generalization of quadratic Gauss sums.  Quadratic Gauss sums are based on primitive $nth$ roots of unity.
The famous formula of Gauss for these sums is the following:
\\
\\
\[G_n = \sum_{k=1}^n exp(2 \pi ik^2/n)\]
\\
\[if~ n\equiv0~(mod~4)~then~G_n = (1+i)\sqrt{n}\]
\[if~ n\equiv1~(mod~4)~then~G_n = \sqrt{n}\]
\[if~ n\equiv2~(mod~4)~then~G_n = 0\]
\[if~ n\equiv3~(mod~4)~then~G_n = i\sqrt{n}\]
\\
\\
However, these sums can be generalized to any finite field.  As is shown in [1], in the case of $GF(p)$ we can define an analogous 
version of these sums which we refer to as quadratic sums in $GF(p)$.  Then the following theorem is true for any prime $p$
and integer $g$ with $gcd(g,p)=1$.
\\
\\
\[Q(g,p) = \sum_{k=1}^n g^{k^2}~(mod~ p)\]
\\
\[where~n~is~ the ~minimum ~positive~ integer~ such~ that\]
\[g^n\equiv1~(mod~p)\]

\[if~ n\equiv0~(mod~4)~then~Q(g,p) \equiv (\sqrt{n}+\sqrt{-n})~(mod~p)\]
\[if~ n\equiv1~(mod~4)~then~Q(g,p) \equiv \sqrt{n}~(mod~ p)\]
\[if~ n\equiv2~(mod~4)~then~Q(g,p) \equiv 0~(mod~ p)\]
\[if~ n\equiv3~(mod~4)~then~Q(g,p) \equiv \sqrt{-n}~(mod~ p)\]
\\

Using the formula for quadratic sums in $GF(p)$ one can calculate square roots of certain integers modulo $p$ in certain cases.
For example, if $n$ is a divisor of $p-1$ and if $n\equiv 1~(mod~4)$ then if given a primitive $nth$ root of unity for the multiplicative 
group of $GF(p)$ one can calculate the square root of $n$ modulo $p$.  The following is an example of this.
\\
\\
\\
\textbf{Example 1}
\\
\\
Suppose that one wishes to calculate the square root of $5$ modulo $41$.  In this case, $18$ generates a subgroup of order $5$
modulo $41$.  Thus using the previously mentioned formula, we have the following result:
\\
\[Q(18,41)\equiv 18^1+18^4+18^9+18^{16}+18^{25}~(mod~41)\]
\[\equiv 18+16+16+18+1\equiv 28~(mod~41)\]

And thus $28$ is a square root of $5$ modulo $41$.
\\
\\
The problem with this method is that adding up the $n$ terms of the quadratic sum one by one is inefficient and should actually be
considered an exponential time algorithm.  In this previous example, it was practical because $n=5$ was a reasonably small integer.
However for most values of $n$ this method would be totally impratical, unless there is a more efficient algorithm for computing
quadratic sums in $GF(p)$.

This brings up an interesting question: Is there a polynomial time algorithm for calculating these types of sums?  If so, it would have
implications for the security of cryptosystems based on the Diffie-Hellman problem in finite fields.  The following section gives a generalization
of the function $Q(g,p)$ and shows how it is closely related to the Diffie-Hellman problem.

\section*{2.1 Quadratic Sums and Diffie-Hellman}
One generalization of the function $Q(g,p)$ is the following which we also refer to as a quadratic sum in $GF(p)$:

\[Q(g,h,p) = \sum_{k=1}^{p-1} g^{k^2}h^k~(mod~ p)\]

One of the most important unsolved problems in number theory or cryptography is to find an efficient algorithm, i.e. a polynomial time
algorithm, for calculating the function $Q(g,h,p)$.  The most obvious way to calculate $Q(g,h,p)$ would be to calculate each of the $p-1$
terms seperately and then add them together.  This would result in an $O(p~log^3p)$ algorithm which is an exponential time algorithm and 
very inefficient.  By polynomial time, we would mean as a polynomial function of $log~p$.

While there is currently no known algorithm for calculating the function $Q(g,h,p)$ in polynomial time, a related case has been solved by Hiary, 
that of calculating truncated theta functions.
\\
\[F_d(a,b,n) = \sum_{k=0}^d exp(2 \pi i (ak^2+bk)/n)\]
\\
In [4] and [5] Hiary gives a polynomial time algorithm for calculating the theta function $F_d(a,b,n)$ as this is useful for calculating the 
Riemann zeta function in certain cases.  In [7] Kuznetzov simplifies Hiary's algorithm using the Mordell integral.  The function $F_d(a,b,n)$ is
quite similar to the function $Q(g,h,p)$.  The main difference is that the first is based on exponentiation involving primitive $nth$ roots of unity
and the second is based on exponentiation in $GF(p)$. 

The function $Q(g,h,p)$ has applications for calculating square roots via the algorithm described in [1].  But more importantly it 
has potential applications concerning the integer factorization problem and the discrete logarithm problem in $GF(p)$.  However the most obvious
application is to the Diffie-Hellman problem in $GF(p)$ which has been conjectured to be equivalent to the discrete logarithm problem.  The 
following explains how the function $Q(g,h,p)$ can be used to solve the Diffie-Hellman problem in $GF(p)$.
\\
\\
\textbf{Theorem 1}
\
\[let~n~be~ the ~minimum ~positive~ integer~ such~ that\]
\[g^n\equiv1~(mod~p)\]
\[if~n\not\equiv 2~(mod~4)~and\]
\[if~h\equiv g^a~(mod~p)~then\]
\[g^{a^2} \equiv Q(g,1,p)(Q(g,h^2,p))^{-1}~(mod~p)\]
\\
The Diffie-Hellman problem or the Computational Diffie-Hellman problem in $GF(p)$ is to calculate the value of $g^{ab}~(mod~p)$ if given
the following three values: $(g,g^a~mod~p,g^b~mod~p)$.  The following formula explains how given that one can calculate  the value of
$g^{a^2}~(mod~p)$ that this can be used to calculate the value of $g^{2ab}~(mod~p)$.
\
\
\[g^{2ab} \equiv (g^{(a+b)^2})(g^{a^2})^{-1}(g^{b^2})^{-1}~(mod~p)\]
\\
The solution to the Diffie-Hellman problem $g^{ab}~(mod~p)$ can be determined by calculating the two square roots of $g^{2ab}~(mod~p)$ 
and then determining which square root represents the correct solution.  What this means is that if one could calculate the quadratic sum $Q(g,h,p)$ 
in polynomial time given any values $g$, $h$, and $p$ then one could also solve the Diffie-Hellman problem in $GF(p)$ in polynomial time.
\\

\section*{3. The Cipolla-Lehmer Square Root Algorithm}

The following explains the algorithm for calculating the function $CL(c,b,p)$.  This definition differs slightly from the algorithm in [9] in that this algorithm
returns a square root of a quadratic residue $c$ in $GF(p)$ or it returns $0$ if the quadratic polynomial selected by the algorithm is not irreducible.
The algorithm in [9] keeps selecting a random quadratic polynomial until an irreducible one is found.
\\
\\
\\
\textbf{Algorithm 1}
\\
\\
\textbf{The Cipolla-Lehmer square root algorithm}
\\
\\
\textbf{Input}: a prime $p$ where $p>2$, a quadratic residue $c$ in $GF(p)$
and an integer $b$ where $0<b<p$
\\
\textbf{Output}: $y$ where $CL(c,b,p) = y$. The output $y$ will be $0$ or a square root of $c$ in $GF(p)$.
\\
\\
(1) $h := (b^2-4c)^{(p-1)/2}~(mod~p)$
\\
(2) if $h \equiv 1~(mod~p)$ or if $h \equiv 0~(mod~p)$ then $s:=0$
\\
(3) if $h \equiv -1~(mod~p)~ then~ s:=1$
\\
(4) $q(x) := x^{(p+1)/2}~ mod~ \langle x^2- bx + c \rangle$ where 
$q(x) = c_1 x + c_0$ for integers $c_0, c_1$
\\
(5) $y:=sc_0$
\\
(6) Return $y$ as the output
\\
\\
\\
\textbf{Example 2}
\\
\\
The following is an example of using Algorithm 1 to calculate $CL(20,2,31)$
\\
\\
\\
(1) $h=(2^2 - (4)(20))^{15}~(mod~31) \equiv 17^{15} \equiv -1~(mod~31)$
\\
(2) Not applicable since $h \equiv -1~(mod~31)$
\\
(3) $s:=1$
\\
(4) $q(x):=x^{16} ~mod~\langle x^2 + 29x+20 \rangle \equiv 0x+19 ~mod~\langle x^2 + 29x+20 \rangle$
\\
(5) $y:=(1)(19)$
\\
(6) Return $19$ as the output
\\
\\
\\
Thus $CL(20,2,31) = 19$ which means that $19$ is a square root of $20$ mod $31$.
\\
\\

\section*{4. The New Square Root Algorithm}

This new algorithm calculates a function $S(d,b,p)$ which calculates the square root of a quadratic residue $d$ in $GF(p)$  based on a random
parameter $b$ or it returns the value of $0$.  Given a fixed quadratic residue $d$ and a fixed prime $p$ where $p \equiv 5~(mod~6)$ consider
the following set of $p-1$ values $for~ 0<k<p : y_k = S(d,k,p)$.  For approximately $1/3$ of these values this algorithm will return $0$.
About $1/3$ of the time it will return $c$ and about $1/3$ of the time it will return $-c~(mod~p)$ where $c$ is the minimum positive integer 
such that $c^2 \equiv d~(mod~p)$.  As an example of this, consider the output of $S(5,k,11)~ for~ 0<k<11$ which calculates a square root of $5$
in $GF(11)$ or returns a value of $0$ based on the parameter $k$.  If $k=1$ or if $k=10$ then $S(5,k,11)=0$.  If $k=2,4,5,$ or $8$ then
$S(5,k,11) = 7$.  If $k=3,6,7,$ or $9$ then $S(5,k,11)=4$.

The algorithm is based on exponentiating modulo a cubic polynomial $f(x)$ in $GF(p^3)$ where $f(x)=x^3+ax+b$ and where the integer values 
$d,b$ and $p$ are given: $p$ is any prime greater than or equal to $5$, $d$ is any nonzero quadratic residue modulo $p$, and $b$ is any integer.
Then the integer $a$ is selected such that $d \equiv -(4a^3+27b^2)~(mod~p)$.  The variable $d$ thus refers to the discriminant of the cubic
polynomial $f(x)$.

If $p \equiv 1~(mod~6)$ then in some cases no such value for the variable $a$ exists, in which case this algorithm will not work.  However,
approximately $1/3$ of the time for randomly selected $d,b$ and $p$ a value of $a$ does exist and so this algorithm will work.  If $p\equiv 5~(mod~6)$
then in all cases a value for $a$ does exist since all integers modulo $p$ are cubic residues.  The following is the main theorem upon which the algorithm
for calculating $S(d,b,p)$ is based.
\\
\\
\\
\\
\textbf{Theorem 2}
\\
\\
\\
Given two integers $a$ and $b$ and a prime $p \ge 5$ such that $gcd(a,p) =1$
\\
Such that the cubic polynomial $x^3+ax+b$ is irreducible modulo $p$
\\
Then the following congruence is true:
\\
\[t^2 \equiv -(4a^3+27b^2)~(mod~p)\]
\[where~ t \equiv (3a)(c_2)^{-1}~(mod~p)\]
\[where~ c_0, c_1~ and~ c_2 ~are~ defined~ as ~any~integers\]
\[such~that~x^p \equiv c_2x^2 + c_1x + c_0~ mod~ \langle x^3 + ax + b \rangle\]
\\
\\
\textbf{Example 3}
\\
\\
\\
The following is an example of using Theorem 2 to calculate the square root of $23$ in $GF(101)$.  Consider the following polynomial:
\[x^3+37x+26\]
This polynomial is irreducible modulo $101$, thus its discriminant $D$ is a quadratic residue.  Using Theorem 2 it follows that:
\\
\[D \equiv t^2 \equiv -((4)(37)^3 + (27)(26)^2)~(mod~101)\]
\[\equiv-(6+72) \equiv 23~(mod~101)\]
\\
by exponentiating in $GF(101^3)$ it follows that:
\\
\[x^{101}\equiv (68x^2+22x+95) ~mod~\langle x^3 + 37x + 26 \rangle \]
thus
\\
\[t \equiv (3a)(c_2)^{-1} \equiv (3)(37)(68)^{-1} \equiv 15~(mod~101) \]
\\
and thus $15$ is a square root of $23$ modulo $101$.
\\
\\
\\
The following theorem is a generalization of Theorem 2 that applies to any irreducible cubic polynomial.
\\
\\
\\
\\
\textbf{Theorem 3}
\\
\\
Given three integers $b,c$ and $d$ and a prime $p \ge 5$ such that $gcd(b^2-3c,p)=1$
and such that the cubic polynomial $x^3+bx^2+cx+d$ is irreducible mod $p$
then the following congruence is true:
\\
\[t^2 \equiv ((18bcd-4b^3d+b^2c^2)-(4c^3+27d^2))~(mod~p)\]
\[where~ t \equiv (b^2-3c)(c_2)^{-1}~(mod~p)\]
\[where~ c_0, c_1~ and~ c_2 ~are~ defined~ as ~any~integers\]
\[such~that~x^p \equiv c_2x^2 + c_1x + c_0~ mod~ \langle x^3 + bx^2 + cx + d \rangle\]
\\
Both Theorem 2 and Theorem 3 show that by exponentiating in $GF(p^3)$, that is, exponentiating modulo irreducible cubic polynomials where the
coefficients are taken modulo some prime $p$, that the square root of the discriminant of the cubic polynomial can be determined.  Both theorems  
define a value for $t$ where $t^2 \equiv D~(mod~p)$ and where $D$ is the cubic polynomial's discriminant.  See [8] for more information on 
the discriminant.
\\
\\
\textbf{Example 4}
\\
\\
The following is an example of using Theorem 3 to calculate the square root of $2$ in $GF(47)$.  Consider the following polynomial:
\[x^3+5x^2+7x+19\]
This polynomial is irreducible modulo 47, thus its discriminant D is a quadratic residue.  Using Theorem 3 it follows that:
\[D \equiv t^2 \equiv ((18)(5)(7)(19) - 4(5)^3(19)+(5)^2(7)^2)-((4)(7)^3+(27)(19)^2)~(mod~47)\]
\[\equiv ((32-6+3)-(9+18)) \equiv 2~(mod~47)\]
by exponentiating in $GF(47^3)$ it follows that:
\[x^{47} \equiv (14x^2+2x+13) ~mod~ \langle x^3 + 5x^2 + 7x + 19 \rangle \]
thus
\[t \equiv (b^2-3c)(c_2)^{-1} \equiv ((5)^2 - 3(7))(14)^{-1} \equiv (4)(37) \equiv 7~(mod~47) \]
and thus $7$ is a square root of $2$ modulo $47$.
\\
\\
Based on Theorem 2, we will next define a function $S(d,b,p)$ that calculates square roots in $GF(p)$ and give an algorithm for calculating it.
\\
\\
\\
\textbf{Definition 1}
\\
\\
\textbf{Definition of $S(d,b,p)$ for $p \equiv 5~(mod~6)$}
\\
\\
Let $p$ be a prime such that $p \equiv 5~(mod~6)$ and let $a$ be the unique solution to the following congruence:
\[d \equiv -(4a^3+27b^2)~(mod~p) \]
If $x^3+ax+b$ is not irreducible modulo $p$ then let $S(d,b,p)=0$.   Otherwise, let $S(d,b,p)=t$ where $t$ is defined in Theorem 2.
\\
\\
\\
\\
\textbf{Definition 2}
\\
\\
\textbf{Definition of $S(d,b,p)$ for $p \equiv 1~(mod~6)$}
\\
\\
Let $p$ be a prime such that $p \equiv 1~(mod~6)$ and let $a$ be any solution to the following congruence (if a solution exists):
\[d \equiv -(4a^3+27b^2)~(mod~p) \]
If no solution $a$ to the above congruence exists or if a solution does exist and
$x^3+ax+b$ is not irreducible modulo $p$ then let $S(d,b,p)=0$.   Otherwise, let $S(d,b,p)=t$ where $t$ is defined in Theorem 2.
\\
\\
\\
It might appear that Definition 2 is ambiguous since $p \equiv 1~(mod~6)$ if a solution $a$ to the previous congruence exists, there will be 
three possible values for $a$ and this definition does not specify which of these three values to use.  However, regardless of which cubic root 
is used for $a$ the same value will be computed for $S(d,b,p)$.

Next we will give an algorithm for calculating the function $S(d,b,p)$ for $p \equiv 5~(mod~6)$ based on Definition 1 and using Theorem 2.  The
most time consuming parts of this algorithm are steps 2, 3 and 4.  Step 2 calculates a cube root in $GF(p)$ and step 3 exponentiates in $GF(p^3)$
and step 4 calculates a multiplicative inverse in $GF(p)$.  All three of these steps each take $O(M(p)~log~p)$ time to calculate.  Thus the whole
algorithm runs in $O(M(p)~log~p)$ time.

The algorithm could be modified to work in the case that $p \equiv 1~(mod~6)$ by checking if $j$ from step 1 is a cubic residue.  If $j$ is a cubic
nonresidue, the algorithm should return $0$, otherwise step 2 should be replaced with an algorithm for calculating cube roots in $GF(p)$ for
$p \equiv 1~(mod~6)$ such as the algorithm in [9].  The rest of the algorithm would remain the same.
\\
\\
\\
\\
\textbf{Algorithm 2}
\\
\\
\textbf{The $GF(p^3)$ square root algorithm for $p \equiv 5~(mod~6)$}
\\
\\
\textbf{Input}: a prime $p$ where $p \equiv 5~(mod~6)$, a quadratic residue $d$ in $GF(p)$
and an integer $b$ where $0<b<p$
\\
\textbf{Output}: $t$ where $S(d,b,p) = t$. The output $t$ will be $0$ or a square root of $d$ in $GF(p)$.
\\
\\
(1) $j :=(d+27b^2)(-4)^{-1}~(mod~p)$
\\
(2) $a:=j^{(2p-1)/3}~(mod~p)$
\\
(3) $q(x):= x^p~mod~ \langle x^3 + ax + b \rangle$ where $q(x) = c_2x^2+c_1x+c_0$
\\
$~~~~~$ for some integers $c_0,~c_1$ and $c_2$
\\
(4) If $x^3+ax+b$ is irreducible in $GF(p)$ then $t:=(3a)(c_2)^{-1}~(mod~p)$
\\
(5) If $x^3+ax+b$ is not irreducible in $GF(p)$ then $t:=0$
\\
(6) Return $t$ as the output
\\
\\
\textbf{Example 5}
\\
\\
The following is a specific example of using Algorithm 2 to calculate $S(21,10,41)$
\\
\\
(1) $j := ((21+27(10)^2)(-4)^{-1}~(mod~41) \equiv(21+27(18))(10) \equiv27~(mod~41)$
\\
(2) $a:=27^{27} \equiv3~(mod~41)$
\\
(3) $q(x) := x^{41}~mod~ \langle x^3 + 3x + 10 \rangle \equiv 30x^2 + 34x + 19 ~mod~ \langle x^3+3x+10 \rangle $
\\
(4) Since $x^3+3x+10$ is irreducible in $GF(41)$ 
\\
$~~~~~$ let $t=(3)(3)(30)^{-1} \equiv 29~(mod~41)$
\\
(5) Not applicable
\\
(6) Return 29 as the output
\\
\\
thus $S(21,10,41)=29$ which means that $29$ is a square root of $21$ mod $41$.

\section*{5. Conjectures involving the function S(d,b,p) }

The following are four conjectures concerning the function $S(d,b,p)$.  Based on calculations that have been done with an implementation of 
Algorithm 2 written in C, it seems probable that these conjectures are true in most if not all cases.
\\
\\
\textbf{Conjecture 1}
\\
\\
If $p$ is a prime such that $p \equiv 5~(mod~6)$ and if $d_1$ is a nonzero quadratic residue modulo $p$ and if
\[d_2 \equiv 729(d_1)^{-1}~(mod~p)\]
then
\\
$~~~~~~~~$ (a) if $S(d_1,1,p)=0$ then $S(d_2,1,p)=0$
\\
$~~~~~~~~$  (b) if $ S(d_1,1,p) \ne 0$ then $S(d_1,1,p)S(d_2,1,p) \equiv -27~(mod~p)$
\\
\\
\textbf{Conjecture 2}
\\
\\
If $p$ is a prime such that $p \equiv 5~(mod~6)$ and $d_2 \equiv b^2d_1(mod~p)$ where $b$ is any integer and $d_1$ and $d_2$ are 
nonzero quadratic residues in $GF(p)$ then
\\
\[(b)(S(d_1,1,p)) \equiv S(d_2,b,p)~(mod~p)\]
\\
\\
\textbf{Conjecture 3}
\\
\\
\\
If $p$ is a prime such that $p \equiv 5~(mod~6)$ then
\\
\\
$~~~~~~~~$ (a) if $p \equiv 2~(mod~9)$ then $S(9,1,p) \equiv -3~(mod~p)$ 
\\
$~~~~~~~~~~~~~$ and $S(81,1,p) \equiv 9~(mod~p)$
\\
\\
$~~~~~~~~$ (b) if $p \equiv 5~(mod~9)$ then $S(9,1,p) \equiv 3~(mod~p)$ 
\\
$~~~~~~~~~~~~~$ and $S(81,1,p) \equiv -9~(mod~p)$
\\
\\
$~~~~~~~~$ (c) if $p \equiv 8~(mod~9)$ then $S(9,1,p) = 0$ 
\\
$~~~~~~~~~~~~~$ and $S(81,1,p) =0$
\\
\\
\\
The following theorem is due to L.E. Dickson [3] (also see [2] and [9]) and gives criteria that can be used to determine whether or not a cubic
polynomial of the form $x^3+ax+b$ is irreducible in $GF(p)$.
\\
\\
\\
\textbf{Theorem 4}
\\
\\
\\
If $p$ is a prime $\ge 5$ and if $f(x) = x^3+ax+b$  for any integers $a$ and $b$ then $f(x)$ is irreducible in $GF(p)$ if and only if the following
two conditions are true:
\\
\\
$~~~~~~$ (a) $D$ is a nonzero quadratic residue in $GF(p)$ where $D = -(4a^3+27b^2)$
\\
\\
$~~~~~~$ (b) $(d_1x+d_2)^{(p^{2}-1)/3} \not\equiv 1~mod~\langle x^2 + 3 \rangle$
\\
$~~~~~~~~~~~~$ where $d_1 \equiv 18^{-1}t~(mod~p)$ and $t^2 \equiv D \equiv -(4a^3 + 27b^2)~(mod~p)$
\\
$~~~~~~~~~~~~$ and $d_2 \equiv -2^{-1}b~(mod~p)$
\\
\\
The following conjecture shows how the function $S(d,b,p)$ can give more specific information about Theorem 4
\\
\\
\\
\\
\textbf{Conjecture 4}
\\
\\
If $p$ is a prime such that $p \equiv 5~(mod~6)$ and if $f(x) = x^3+ax+b$ for any integers $a$ and $b$ such that $gcd(ab,p)=1$ then the following
two conditions are true:
\\
\\
\\
$~~~~~~~$ (a) $S(-(4a^3 + 27b^2),b,p)=0$ if and only if 
\\
$~~~~~~~~~~~~$ $f(x)$ is not irreducible in $GF(p)$
\\
\\
$~~~~~~~$ (b) if $f(x)$ is irreducible in $GF(p)$ then
\\
$~~~~~~~~~~~~$ $(d_1x+d_2)^{(p^{2}-1)/3} \equiv -2^{-1}(x+1)~mod~\langle x^2 + 3 \rangle$
\\
$~~~~~~~~~~~~$ where $d_1 \equiv 18^{-1}t~(mod~p)$ and $t^2 \equiv D \equiv -(4a^3 + 27b^2)~(mod~p)$
\\
$~~~~~~~~~~~~$ where $t$ is defined as $t=S(D,b,p)$
\\
$~~~~~~~~~~~~$ and $d_2 \equiv -2^{-1}b~(mod~p)$
\\

\section*{5.1 Two Cubic Residuosity Conjectures }

We define the concept of a residuosity theorem as the following: given two primes $p$ and $q$, and a function $f(q,p)$ that is computed modulo $p$,
a residuosity theorem is any theorem that shows a relationship between the output of the function $f(q,p)$ and the value of $p^{(q-1)/c}~(mod~q)$
where $c>1$ and $q \equiv 1~(mod~c)$.  If $c=2$ this would be a quadratic residuosity theorem.  If $c=3$, this would be a cubic residuosity theorem.

Using this definition the most well-known residuosity theorem would be the law of quadratic reciprocity which given two odd primes $p$ and $q$ shows a
relationship between the value of $p^{(q-1)/2}~(mod~q)$ and the value of $q^{(p-1)/2}~(mod~p)$.  In this case, the function $f(q,p)$ would be 
defined as $f(q,p) = q^{(p-1)/2}~(mod~p)$.  In the following section we will present two cubic residuosity conjectures concerning the function $S(d,b,p)$
for $b=1$ and for $b=2$ that seem to be true based on computational evidence.
\\
\\
\textbf{Conjecture 5}
\\
\\
If $p$ is a prime such that $p \equiv 5~(mod~6)$ and if $d \equiv 81e^2~(mod~p)$ for some integer $e$ such that $gcd(e,p)=1$ with 
$d \not\equiv 9~(mod~p) $ and$~d \not\equiv 81~(mod~p) $
\\
And suppose that the following criteria are met for any positive integers $x$ and $y$:
\\
\\
(1) $x \equiv (e-1)(2)^{-1}~(mod~p)$
\\
(2) $y \equiv (e+1)(2)^{-1}~(mod~p)$
\\
(3) $x \equiv 1~(mod~3)$
\\
(4) $y \equiv 2~(mod~3)$
\\
(5) $q$ is a prime where $q = x^2 +xy + y^2$
\\
\\
Then the following three statements are true:
\\
\\
(a) $S(d,1,p) = 0$ if and only if $p^{(q-1)/3} \equiv 1~(mod~q)$
\\
(b) $S(d,1,p) \equiv 9e~(mod~p)$  if and only if $p^{(q-1)/3} \equiv xy^{-1}~(mod~q)$
\\
(c) $S(d,1,p) \equiv -9e~(mod~p)$  if and only if $p^{(q-1)/3} \equiv x^{-1}y~(mod~q)$
\\
\\
\\
\textbf{Example 6}
\\
\\
Suppose that $p$ is a prime $p \equiv 5~(mod~6)$ and that $d \equiv 729~(mod~p)$ 
\\
thus $e^2 \equiv 9~(mod~p)$
\\
\\
Using criteria (1) - (4) in Conjecture 5 we could choose $x=1$ and $y=2$
\\
and thus by (5) $q = 1^2 + (1)(2) + 2^2 = 7$
\\
\\
The output of $S(d,1,p)$ would depend on the value of $p^{(7-1)/3} \equiv p^2~(mod~7)$
\\
\\
Thus items a, b and c from Conjecture 5 would imply the following:
\\
\\
(a) $S(d,1,p) = 0$ if and only if $p \equiv \pm 1~(mod~7)$
\\
(b) $S(d,1,p) \equiv 27~(mod~p)$ if and only if $p \equiv \pm 2~(mod~7)$
\\
(c) $S(d,1,p) \equiv -27~(mod~p)$ if and only if $p \equiv \pm 4~(mod~7)$
\\
\\
\\
\\
\textbf{Conjecture 6}
\\
\\
If $p$ is a prime such that $p \equiv 5~(mod~6)$ and if $d \equiv 81e^2~(mod~p)$ for some integer $e$ such that $gcd(e,p)=1$ with 
$d \not\equiv 9~(mod~p) $ and$~d \not\equiv 81~(mod~p) $
\\
And suppose that the following criteria are met for any positive integers $x$ and $y$:
\\
\\
(1) $x \equiv (e-2)~(mod~p)$
\\
(2) $y \equiv (e+2)~(mod~p)$
\\
(3) $x \equiv 1~(mod~3)$
\\
(4) $y \equiv 2~(mod~3)$
\\
(5) $q$ is a prime where $q = x^2 +xy + y^2$
\\
\\
Then the following three statements are true:
\\
\\
(a) $S(d,2,p) = 0$ if and only if $p^{(q-1)/3} \equiv 1~(mod~q)$
\\
(b) $S(d,2,p) \equiv 9e~(mod~p)$  if and only if $p^{(q-1)/3} \equiv xy^{-1}~(mod~q)$
\\
(c) $S(d,2,p) \equiv -9e~(mod~p)$  if and only if $p^{(q-1)/3} \equiv x^{-1}y~(mod~q)$
\\
\\
\\
\\
\textbf{Example 7}
\\
\\
Suppose that $p$ is a prime $p \equiv 5~(mod~6)$ and that $d \equiv 729~(mod~p)$ 
\\
thus $e^2 \equiv 9~(mod~p)$
\\
\\
Using criteria (1) - (4) in Conjecture 6 we could choose $x=1$ and $y=5$
\\
and thus by (5) $q = 1^2 + (1)(5) + 5^2 = 31$
\\
\\
The output of $S(d,2,p)$ would depend on the value of $p^{(31-1)/3} \equiv p^{10}~(mod~31)$
\\
\\
The following are three specific examples of what items a, b and c from Conjecture 6 would imply:
\\
\\
(a) $S(d,2,23) = 0$ since $23^{10} \equiv  1~(mod~31)$
\\
(b) $S(d,2,59) \equiv 27~(mod~59)$ since $59^{10} \equiv (1)(5)^{-1} \equiv 25~(mod~31)$
\\
(c) $S(d,2,41) \equiv -27~(mod~41)$ since $41^{10} \equiv (1)^{-1}(5) \equiv 5~(mod~31)$
\\
\\
\\
One may notice that both Conjectures 5 and 6 are very similar to each other.  Conjecture 5 pertains to the case $b=1$ and Conjecture 6 to the 
case $b=2$.  This could be extended to the cases of $b=p-1$ and $b=p-2$ if one notes the identity: $S(d,b,p) \equiv -S(d,p-b,p)(mod~p)$.
One might suspect that it could be possible to give a generalization for other values of $b$.  In fact such a generalization does exist 
which will be mentioned in the next section 5.2.  This generalization can be considered to be a cubic reciprocity formula.

Also one should note that both Examples 6 and 7 considered the simplest case of $d \equiv 729~(mod~p)$ which allowed for the smallest possible 
value of $q$ which was $7$ for Conjecture 5 and $31$ for Conjecture 6.
\\

\section*{5.2 Cubic Reciprocity}

Both conjectures 5 and 6 can be generalized and these generalizations can be considered to be cubic reciprocity identities.  By cubic reciprocity identity we would mean that if given two primes $p$ and $q$ the problem of determining whether or not $q$ is a cubic residue in $GF(p)$ can be shown to be equivalent to making some determination about an element of a group that uses arithmetic modulo $q$ instead of arithmetic modulo $p$.  One of the oldest examples is Jacobi's rational cubic reciprocity theorem which is mentioned in [12].  What is meant by rational cubic reciprocity theorem is that it can be used to determine whether or not $q$ is a cubic residue mod $p$.  However in the case that $q$ is a cubic nonresidue in $GF(p)$ the theorem gives no indication as to which primitve cubic root of unity the value $q^{(p-1)/3}~(mod~p)$ is equal to.  The following is a slightly changed version of Jacobi's rational cubic reciprocity theorem as presented in [12].
\\
\\
\textbf{Theorem 5}
\\
\\
If $p$ and $q$ are primes where $p \equiv 1~ (mod~ 3)$ and $q \equiv 1~ (mod~ 3)$, $q \ne p$ and $4p = L^{2}+27M^{2}$  and $4q = A^{2}+27B^{2}$ then $q$ is a cubic residue in $GF(p)$ if and only if $(LB-AM)(LB+AM)^{-1}$  is cubic residue in $GF(q)$
\\
\\
Theorem 5 can be generalized so that it is a full cubic reciprocity theorem instead of just a rational cubic reciprocity identity.  The following is an example of this from [12].
\\
\\
\textbf{Theorem 6}
\\
\\
If $p$ and $q$ are primes where $p \equiv 1~ (mod~ 3)$ and $q \equiv 1~ (mod~ 3)$, $q \ne p$ and $4p = L^{2}+27M^{2}$  and $4q = A^{2}+27B^{2}$ where $A$ is defined so that $A \equiv 1~(mod~3)$ then $q^{(p-1)/3} \equiv ((-1-L/(3M))2^{-1})^{i}~(mod~p)$ if and only if
\\
$((LB-AM)(LB+AM)^{-1})^{(q-1)/3}$   $\equiv ((-1-A/(3B))2^{-1})^{i}~(mod~q)$ 
\\
\\
Both Theorem 5 and it's generalization Theorem 6 represent very efficient cubic reciprocity identities.  Determining whether or not $q$ is a cubic residue $mod~ p$ which would normally require $O(log^{3}~p)$ time to calculate can instead be translated into determining whether or not a certain element is a cubic residue $mod~ q$ which runs in $O(log^{3}~q)$ time.  If $q$ much smaller than $p$ then this can be significantly more efficient.  The problem is that both theorems require that $q \equiv 1~(mod~3)$ and so cannot be used if the prime $q \equiv 2~(mod~3)$.

There are some cubic reciprocity theorems that can be used in the case that $q \equiv 2~(mod~3)$ such as the rational cubic reciprocity theorem due to Emma Lehmer which is stated as Theorem 1.2 in [13] or the similar but more general full cubic reciprocity Theorem 1.2 in [12].  Another very interesting example is given by Zhi-Hong Sun in [12] which is based on
exponentiation in a unique type of group called $C(q)$ which is  defined on the set of all integers modulo $q$ with the exception of the two square roots of $q-3$ if $q \equiv 1~(mod~3)$ plus an identity element which we denote by $e$.  The group operation $*$ where $x$ and $y$ refer to any two elements of the group both not equal to $e$ is defined as the following:
\\
\\
(1) $e*e = e$
\\
(2) $e*x = x*e = x$
\\
(3) if $x \equiv -y~(mod~q)$ then $x*y = e$
\\
(4) if $x \not\equiv -y~ (mod~q)$ then $x*y = (xy-3)(x+y)^{-1}~(mod~q)$
\\
\\
Using the definition of $C(q)$ we have the following Theorem from [12]. 
\\
\\
\textbf{Theorem 7}
\\
\\
Let $p$ and $q$ be distinct primes both greater than 3 and suppose $p \equiv ~1~(mod~3)$ and that $4p = L^{2}+ 27M^{2}$ then
$q$ is a cubic residue if $GF(p)$ if and only if
$c$ is a cube in $C(q)$ where $c \equiv L(3M)^{-1}~ (mod~q)$ 

The order of the group $C(q)$ where $q$ is prime greater than $3$ is $q-1$ if $q \equiv 1~(mod~3)$ and is $q+1$ if $q \equiv 2~(mod~3)$.  Thus determining whether or not an element is a cubic residue involves exponentiating to the $(q-1)/3$ or $(q+1)/3$ power depending on the value of $q ~mod~3$.  One should notice that the group operation $*$ involves calculating a multiplicative inverse in $GF(p)$ which takes $O(log^{3}~q)$ time.  Thus exponentiating in the group $C(q)$ takes $O(log^{4}~q)$ time.

In all the previously  mentioned cubic reciprocity formulas, it was required to know solutions $L$ and $M$ to the equation $4p = L^{2}+27M^{2}$.  In the next cubic reciprocity formula based on the $GF(p^{3})$ square root function $S(d,b,p)$ the formula requires knowing the values of $x$ and $y$ where $x^{2}+xy+y^{2} = p$.  However these two equations are essentially equivalent.  Any prime $p$ that is congruent to $1~(mod~3)$ can be written in the form  $x^{2}+xy+y^{2}$.  Without loss of generality we can assume that $x$ and $y$ are both positive integers and that there are only three cases to consider: (1) $x \equiv 1~(mod~3)$ and $y \equiv 2~(mod~3)$ ~~ (2) $x \equiv 0~(mod~3)$ and $y \equiv 1~(mod~3)$ or (3) $x \equiv 0~(mod~3)$ and $y \equiv 2~(mod~3)$.  If case 1 is true then $L = |y-x|$ and $M = (x+y)/3$.  If cases 2 or 3 are true then $L = x+2y$ and $M = x/3$.

Both conjectures 5 and 6 can be considered to be cubic reciprocity identitiies but only apply to certain cases.  For example, conjecture 5 corresponds to the case $y-x = 1$ where $p = x^{2}+xy+y^{2}$ and where $x \equiv 1~(mod~3)$ and $y \equiv 2~(mod~3)$.  The following is an example of how this can be generalized to be a cubic reciprocity identity for all possible values of $x$ and $y$.
\\
\\
\textbf{Conjecture 7}
\\
\\
Suppose that $p$ is a prime where $p \equiv 1~(mod~6)$ and that $p = x^{2}+xy+y^{2}$ where $x$ and $y$ are both positive integers and suppose $d \equiv 81(x+y)^{2} ~(mod~q)$ and $b \equiv y-x ~(mod~q)$ and that $x^{2} \not\equiv y^{2}~ (mod~q)$ where $q$ is any prime such that $q \equiv 5~(mod~6)$ then
\\
\\  Case 1: $x \equiv 1~(mod~3)$ and $ y \equiv 2~(mod~3)$
\\
\\
(a) $S(d,b,q) = 0$ if and only if $q^{(p-1)/3} \equiv 1~(mod~p)$
\\
(b) $S(d,b,q) \equiv 9(x+y)~(mod~q)$  if and only if $q^{(p-1)/3} \equiv xy^{-1}~(mod~p)$
\\
(c) $S(d,b,q) \equiv -9(x+y)~(mod~q)$  if and only if $q^{(p-1)/3} \equiv x^{-1}y~(mod~p)$
\\
\\ Case 2: $x \equiv 0~(mod~3)$ and $ y \equiv \pm 1~(mod~3)$
\\
\\
Let $t = ((q+1)/3) ~mod~3$
\\
\\
If $S(d,b,q) = 0$ then $c=0$
\\
If $S(d,b,q) \equiv 9(x+y)~(mod~q)$ then $c=1$
\\
If $S(d,b,q) \equiv -9(x+y)~(mod~q)$ then $c=2$   
\\
\\
Let $z = (t+c)~mod~3$
\\
\\
(a) $z = 0$ if and only if $q^{(p-1)/3} \equiv 1~(mod~p)$
\\
(b) $z = 1$  if and only if $q^{(p-1)/3} \equiv xy^{-1}~(mod~p)$
\\
(c) $z = 2$  if and only if $q^{(p-1)/3} \equiv x^{-1}y~(mod~p)$
\\
\\
Conjecture 7 only covers the case $q \equiv 5~(mod~6)$.  It can be generalized to cover the other case $q \equiv 1~(mod~6)$.  However this generalization will not be mentioned here since we did not give a specific method for calculating $S(d,b,q)$ if  $q \equiv 1~(mod~6)$. One should note that the same exact algorithm applies as for the case  $q \equiv 5~(mod~6)$. except that a different algorithm for calculating cube roots has to be applied.  Also one should note that the Conjecture 7 algorithm gives essentially the same output as the method based on the group $C(q)$ but is somewhat faster running in $O(log^{3}~ q)$ time as opposed to $O(log^{4}~q)$ time.

\section*{6. Conclusion}

We presented two previously known methods for calculating square roots in $GF(p)$: the quadratic sum method and the Cipolla-Lehmer method.
Also we showed how quadratic sums are related to the Diffie-Hellman problem in $GF(p)$ and how efficient methods for calculating the function 
$Q(g,h,p)$ might lead to efficient methods for solving the Diffie-Hellman problem in $GF(p)$.

We also presented a new method for calculating square roots in $GF(p)$ based on exponentiation in $GF(p^3)$.  A function $S(d,b,p)$ was defined and
an algorithm for calculating this function was given.  In addition to this, seven conjectures relating to the output of the square root function $S(d,b,p)$
were given.  The seventh conjecture mentioned was a very efficient new cubic reciprocity identity.
This new $GF(p^3)$ square root algorithm like the Cipolla-Lehmer algorithm was shown to run in $O(M(p)~log~p)$ time where $M(p)$ is
the amount of time it takes to calculate one multiplication modulo $p$.

Python 3 and C implementations [14] of the $GF(p^{3})$ square root algorithm are available at https://github.com/davidsknight/NumberTheory.  The program modsqrt.py contains both the $GF(p^{3})$ algorithm as well as the Cipolla-Lehmer algorithm.  The program modsqrt2.c is a C implementation of these two algorithms.  The Python progam cubicreciprocity.py checks the correctness of conjecture 7.  And cubicsqrt.py is a polynomial time deterministic algorithm based on Theorem 3 for determining the square root in $GF(p)$ of the discriminant of a given cubic polynomial that is irreducible $mod~ p$. 
\\
\\
\\
\\
\section*{References }

[1] Eric Bach and Klaus Huber. \textit{Note on taking square-roots modulo $N$}.  IEEE Transactions on Information Theory 45(2):807-809, 1999.
\\
\\
\
[2] Gook Hwa Cho, Namhun Koo, Eunhye Ha, and Soonhak Kwon. \textit{New cube root algorithm based on third order linear recurrence relation in
finite field}, preprint available from http://eprint.iacr.org/2013/024.pdf, 2013.
\\
\\
\
[3] L.E. Dickson. \textit{Criteria for irreducibility of functions in a finite field}.  Bulletin of the American Mathematical Society, Volume 13(1):1-8, 1906.
\\
\\
\
[4] G.A. Hiary. \textit{Fast methods to compute the Reimann zeta function}.  Annals of Mathematics 174(2):891-946, 2011.
\\
\\
\
[5] G.A. Hiary. \textit{A nearly-optimal method to compute the truncated theta function, its derivatives, and integrals}.  Annals of Mathematics 
174(2):859-889, 2011.
\\
\\
\
[6] Namhun Koo, Gook Hwa Cho, and Soonhak Kwon. \textit{Square root algorithm in $F_q$ for $q \equiv 2^s+1~(mod~2^{s+1})$}, preprint available
from 
\\
 http://eprint.iacr.org/2013/087.pdf, 2013.
\\
\\
\
[7]  A. Kuznetsov. \textit{Computing the truncated theta function via Mordell integral}, preprint available from arXiv:math.NT/1306.4081, June 2013.
\\
\\
\
[8] Yang Min, Meng Qingshu, Wang Zhangyi, Li Li, Zhang Huanguo. \textit{Some observations to speed the polynomial selection in the number field
sieve}, preprint available at http://eprint.iacr.org/2012/599.pdf, 2012.
\\
\\
\
[9] Nozomu Nishihara, Ryuichi Harasawa, Yutaka Sueyoshi, Aichi Kudo. \textit{A remark on the computation of cube roots in finite fields}, preprint
available at http://eprint.iacr.org/2009/457.pdf, 2009.
\\
\\
\
[10] Rene Schoof. \textit{Elliptic curves over finite fields and the computation of square roots mod p}.  Mathematics of Computation, 44(170):483-494,
April 1985.
\\
\\
\
[11] Tsz-Wo Sze. \textit{On taking square roots without quadratic nonresidues over finite fields}.  Mathematics of Computation, 802(275):1797-1811, 
July 2011.
\\
\\
\
[12] Zhi-Hong Sun. \textit{Cubic residues and binary quadratic forms}. Journal of Number Theory 124(1):62-104, May 2007.
\\
\\
\
[13]  Zhi-Hong Sun. \textit{On the theory of cubic residues and nonresidues}. Acta Arithmetica 84(4):291-335, Jan 1998.
\\
\\
\
[14] David S. Knight. https://github.com/davidsknight/NumberTheory

\end{document}